# NON-UNIFORM IN TIME STATE ESTIMATION OF DYNAMICAL SYSTEMS

**Iasson Karafyllis**[*] **and Costas Kravaris**[**]

[*] **Dept. of Environmental Engineering, Technical University of Crete**
Email:ikarafyl@enveng.tuc.gr

[**] **Dept. of Chemical Engineering, University of Patras**
Email: kravaris@chemeng.upatras.gr

**Abstract**
In this paper it is showed that if a time-varying uncertain system is robustly completely detectable then there exists an estimator for this system, i.e. we can estimate asymptotically the state vector of the system. Moreover, if a time-varying uncertain system is robustly completely observable then there exists an estimator for this system that guarantees convergence of the estimates with assignable rate of convergence. Finally, it is proved that under the assumption of Robust Lipschitz complete observability, there is a global solution of the observer problem.

## 1. Introduction

One of the biggest challenges of Mathematical Control Theory has been the problem of constructing state observers or estimators for nonlinear systems. This problem has attracted a lot of attention in the literature in the past decades (see for example, [1,2,4,6,8-10,15,17,19-21,25,29,30,32]).

It should be noticed that the problem of design of observers for nonlinear systems has been approached from different research directions. Tsinias in [29,30] provided sufficient conditions for the design of nonlinear global time-invariant observers based on Lyapunov-like characterizations of observability and detectability, which can be verified easily for a special class of nonlinear systems. The works of Gauthier, Kupka and others ([8-10]) have provided semi-global solutions to the observer problem for systems with analytic dynamics based on a high gain strategy. The case of observable systems with analytic dynamics and the solvability of a series solution methodology to the observer problem with assignable exponential rate of convergence in transformed coordinates has been considered initially in [17] and later in [20-21]. A transformed coordinates approach for a limited class of systems with smooth dynamics has provided local solutions to the observer problem in [19,25]. On the other hand, a global solution to the observer problem is provided in [1] under the assumptions of Output-to-State Stability and Global Complete Observability. Observers with delays have been considered in [23] for special classes of nonlinear systems and time-varying observers for linear time-varying systems, which guarantee non-uniform in time convergence, have been considered in [32].

The present paper will provide new results regarding the nonlinear state estimation problem, referring to a broad class of systems (time-varying uncertain nonlinear systems, which of course include autonomous systems as a special case), under minimal regularity conditions (local Lipschitz continuity for the dynamics and simple continuity for the output map) and easily verifiable observability assumptions, leading to global solutions to the observer problem with assignable rate of convergence of the error.

In particular, the present work focuses on the state estimation problem for time-varying nonlinear forward complete systems of the form:

$$\dot{x}(t) = f(t, x(t)) \quad ; \quad y(t) = h(t, x(t))$$
$$x \in \Re^n, t \geq 0, y \in \Re^k \tag{1.1}$$

the main result being that, under the assumption of Lipschitz complete observability for (1.1), there exists a global solution to the observer problem for (1.1). Roughly speaking, this means that if (1.1) is Lipschitz completely observable, then for every smooth non-decreasing function $\phi: \Re^+ \to [1,+\infty)$, there exist **time-varying** vector fields $k: \Re^+ \times \Re^m \times \Re^k \to \Re^m$, $\Psi: \Re^+ \times \Re^m \times \Re^k \to \Re^n$ such that the solution of system (1.1) with

$$\dot{z}(t) = k(t, z(t), y(t)) \quad ; \quad \bar{x}(t) = \Psi(t, z(t), y(t))$$
$$z \in \Re^m, t \geq 0, \bar{x} \in \Re^n \tag{1.2}$$

and initial condition $(x(t_0), z(t_0)) = (x_0, z_0) \in \Re^n \times \Re^m$ satisfies the following **global convergence** property:

$$\lim_{t \to +\infty} \phi(t) |x(t) - \bar{x}(t)| = 0, \quad \forall (t_0, x_0, z_0) \in \Re^+ \times \Re^n \times \Re^m \tag{1.3}$$

and moreover satisfies the **Consistent Initialization Property**, i.e., for every $(t_0, x_0) \in \Re^+ \times \Re^n$ there exists $z_0 \in \Re^m$ such that the solution $(x(\cdot), z(\cdot))$ of system (1.1) with (1.2) initiated from $(x_0, z_0) \in \Re^n \times \Re^m$ at time $t_0 \geq 0$, satisfies

$$x(t) = \bar{x}(t), \quad \forall t \geq t_0$$

The results to be presented in this paper also cover the case of uncertain forward complete nonlinear time-varying systems:

$$\dot{x}(t) = f(t, x(t), d(t)) \quad ; \quad y(t) = h(t, x(t))$$
$$x \in \Re^n, t \geq 0, d(t) \in D, y \in \Re^k \tag{1.4}$$

where $D \subset \Re^l$ is a known compact set. Under appropriate robust complete observability assumptions for (1.4), the global convergence property (1.3) can be imposed, leading to the notion of $\phi$-estimator for (1.4). Roughly speaking, our result states that if system (1.4) is robustly completely observable, then for every smooth non-decreasing function $\phi : \Re^+ \to [1, +\infty)$ there exist time-varying vector fields $k : \Re^+ \times \Re^m \times \Re^k \to \Re^m$, $\Psi : \Re^+ \times \Re^m \times \Re^k \to \Re^n$ such that the solution of system (1.4) with (1.2) satisfies the global property (1.3) for all measurable functions $d : \Re^+ \to D$. It should be emphasized that a $\phi$-estimator in the above sense is **not** an observer, since it does **not** necessarily satisfy the property of **consistent initialization** (see [15,29,30]). On the other hand, the property of consistent initialization cannot in general be satisfied if the original system is uncertain, i.e. its dynamics contain unknown parameters. However, the notion of the estimator is applicable even in this case. Moreover, it is known that dynamic output stabilization methods are based on estimates of the state vector of the system and for such purposes the property of consistent initialization is not essential (see [7,28]). A preliminary version of this paper concerning systems without uncertainties was given in [16].

It must be emphasized that the notions of complete detectability and complete observability (Robust Complete Observability/Detectability) that will be used in this work generalize the corresponding notion of Uniform Complete Observability presented in [28] for autonomous systems, as well as similar notions given in [10]. In particular, for disturbance-free systems with analytic output maps and dynamics, the notion Robust Complete Observability used in the present work coincides with the notion of Uniform Complete/Infinitesimal Observability of [10], for which appropriate test conditions are available.

The contents of the paper are presented as follows. In Section 2 we provide all necessary definitions of the notions used as well as all the technical results needed for the proofs of the main results of the paper. In Section 3 of the present paper we present the main results of the paper: if system (1.4) is robustly completely detectable then there exists an estimator for system (1.4) (Theorem 3.1). Moreover, if system (1.4) is Robustly Lipschitz completely observable, then the observer problem for (1.4) can be solved globally (Theorem 3.3).

All the results given in the present paper are **existence results, proved in a constructive way**, including conditions that guarantee **robustness to modeling errors.** The question of construction of nonlinear observers, which are robust to both measurement noise and modeling errors, as well as implementation issues of the proposed observers/estimators, will be the subject of future research.

**Notation**

* By $M_D$ we denote the set of all measurable functions from $\Re^+$ to $D$, where $D \subset \Re^m$ is a given compact set.
* By $C^j(A)$ ($C^j(A;\Omega)$), where $j \geq 0$ is a non-negative integer, we denote the class of functions (taking values in $\Omega$) that have continuous derivatives of order $j$ on $A$. $L^\infty(A;B)$ ($L^\infty_{loc}(A;B)$) denotes the set of all measurable functions $u : A \to B$ that are (locally) essentially bounded on $A$.

* For $x \in \Re^n$, $x'$ denotes its transpose and $|x|$ its usual Euclidean norm.
* By $B[x,r]$ where $x \in \Re^n$ and $r \geq 0$, we denote the closed sphere in $\Re^n$ of radius $r$, centered at $x \in \Re^n$.
* $x(t) = x(t, t_0, x_0; d)$ denotes the unique solution of (1.4) at time $t \geq t_0$ that corresponds to some input $d(\cdot) \in M_D$, initiated from $x_0 \in \Re^n$ at time $t_0 \geq 0$.
* For the definition of the class $K_\infty$, see [18]. By $KL$ we denote the set of all continuous functions $\sigma = \sigma(s,t): \Re^+ \times \Re^+ \to \Re^+$ with the properties: (i) for each $t \geq 0$ the mapping $\sigma(\cdot, t)$ is of class $K$; (ii) for each $s \geq 0$, the mapping $\sigma(s, \cdot)$ is non-increasing with $\lim_{t \to +\infty} \sigma(s,t) = 0$.
* The saturation function is defined on $\Re$ as $sat(x) := \begin{cases} x & if \quad |x| < 1 \\ x/|x| & if \quad |x| \geq 1 \end{cases}$.

## 2. Definitions and Preliminary Technical Results

In this section we provide definitions and technical lemmas that play a key role in the proofs of the main results of the paper. Their proofs can be found at the Appendix, unless otherwise stated.

**Definition 2.1** *We denote by $K^+$ the class of $C^0$ functions $\phi: \Re^+ \to \Re$ and we denote by $K^* \subset K^+$ the class of non-decreasing $C^\infty$ functions with $\phi(0) \geq 1$, which belong to $K^+$ and satisfy $\lim_{t \to +\infty} \dot{\phi}(t) \phi^{-2}(t) = 0$.*

For example the functions $\phi(t) = 1$, $\phi(t) = 1+t$, $\phi(t) = \exp(t)$ all belong to the class $K^*$. The proof of Lemma 2.2 in [12] actually shows an important property for this class of functions: for every function $\phi$ of class $K^+$, there exists a function $\tilde{\phi}$ of class $K^*$, such that: $\phi(t) \leq \tilde{\phi}(t)$ for all $t \geq 0$. We next give the notion of Robust Forward Completeness, which was introduced in [14] for uncertain dynamical systems. Consider the system (1.4), where $D \subset \Re^l$ is a compact subset and the vector fields $f: \Re^+ \times \Re^n \times D \to \Re^n$, $h: \Re^+ \times \Re^n \to \Re^k$ with $f(t, 0, d) = 0$, $h(t, 0) = 0$ for all $(t, d) \in \Re^+ \times D$, satisfy the following conditions:

1) The functions $f(t, x, d)$, $h(t, x)$ are continuous.
2) The function $f(t, x, d)$ is locally Lipschitz with respect to $x$, uniformly in $d \in D$, in the sense that for every bounded interval $I \subset \Re^+$ and for every compact subset $S$ of $\Re^n$, there exists a constant $L \geq 0$ such that:
$$|f(t, x, d) - f(t, y, d)| \leq L|x - y|$$
$$\forall t \in I, \forall (x; y) \in S \times S, \forall d \in D$$

Let us denote by $x(t, t_0, x_0; d) = x(t)$ the unique solution of (1.4) at time $t$ that corresponds to input $d \in M_D$, with initial condition $x(t_0) = x_0$ and let $h(t, x(t, t_0, x_0; d)) = y(t)$.

**Definition 2.2** *We say that (1.4) is **Robustly Forward Complete (RFC)** if for every $T \geq 0$, $r \geq 0$ it holds that:*
$$\sup\{|x(t_0 + s)|\,;\, |x_0| \leq r,\, t_0 \in [0,T],\, s \in [0,T],\, d(\cdot) \in M_D\} < +\infty$$

The following proposition clarifies the consequences of the notion of Robust Forward Completeness and provides estimates of the solutions. Its proof can be found in [14].

**Proposition 2.3 (Lemma 2.3 in [14])** *Consider system (1.4) with $d \in D$ as input. System (1.4) is RFC if and only if there exist functions $\mu \in K^+$, $a \in K_\infty$ such that for every input $d(\cdot) \in M_D$ and for every $(t_0, x_0) \in \Re^+ \times \Re^n$, the unique solution $x(t)$ of (1.4) corresponding to $d(\cdot)$ and initiated from $x_0$ at time $t_0$ exists for all $t \geq t_0$ and satisfies:*
$$|x(t)| \leq \mu(t)\, a(|x_0|),\ \forall t \geq t_0 \tag{2.1}$$

The notions of robust complete observability and robust complete detectability for time-varying systems are given next. The definitions given here directly extend the corresponding notions given in [28], concerning autonomous systems, as well as similar notions given in [10] for autonomous systems with analytic dynamics.

**Definition 2.4** *Consider the system (1.4) with* $h(t,x) = (h_1(t,x),...,h_k(t,x))$, $h_j \in C^0(\Re^+ \times \Re^n ; \Re)$ *(* $j = 1,...,k$ *) and* $h_j(t,0) = 0$ *for all* $t \geq 0$ *(* $j = 1,...,k$ *). Suppose that (1.4) is RFC. Let* $g_j : \Re^+ \times \Re^k \to \Re$, $a_j : \Re^+ \times \Re^k \to \Re$, $j = 1,...,p$ *be functions of class* $C^0(\Re^+ \times \Re^k ; \Re)$ *with* $\inf\{a_j(t,y) ; (t,y) \in \Re^+ \times \Re^k\} > 0$, $g_j(t,0) = 0$ *for all* $t \geq 0$. *Let* $m_j \geq 0$ *(* $j = 1,...,p$ *) integers with the property that the family of functions, defined recursively below for each* $j = 1,...,p$ *:*

$$y_{0,j}(t,x) = g_j(t,h(t,x))$$

$$y_{i,j}(t,x) := \frac{1}{a_j(t,h(t,x))}\left\{\frac{\partial y_{(i-1),j}}{\partial t}(t,x) + \frac{\partial y_{(i-1),j}}{\partial x}(t,x) f(t,x,d) + \varphi_{i,j}(t,h(t,x))\right\}, \; i = 1,...,m_j$$

*are all independent of* $d \in D$ *and of class* $C^1(\Re^+ \times \Re^n ; \Re)$, *where* $\{\varphi_{i,j}\}$ $i = 1,...,m_j$, $j = 1,...,p$ *are functions of class* $C^0(\Re^+ \times \Re^k ; \Re)$ *such that each* $\varphi_{i,j}(t,h(t,x))$ *is locally Lipschitz with respect to* $x$, *with* $\varphi_{i,j}(t,0) = 0$ *for all* $t \geq 0$. *Let* $\mathcal{D}_j y(t,x) := (y_{1,j},..., y_{m_j,j}) \in \Re^{m_j}$, $j = 1,...,p$ *and* $\mathcal{D} y(t,x) := (\mathcal{D}_1 y(t,x),...,\mathcal{D}_p y(t,x)) \in \Re^m$ *where* $m := \sum_{j=1}^{p} m_j$. *We say that a function* $\theta \in C^0(\Re^+ \times \Re^n ; \Re^l)$ *with* $\theta(\cdot, 0) = 0$ *is **robustly completely observable with respect to (1.4)** if there exists a function* $\Psi \in C^0(\Re^+ \times \Re^k \times \Re^m ; \Re^l)$ *with* $\Psi(t,0,0) = 0$ *for all* $t \geq 0$ *such that*

$$\theta(t,x) = \Psi(t, h(t,x), \mathcal{D} y(t,x)), \; \forall (t,x) \in \Re^+ \times \Re^n \quad (2.2)$$

*We say that system (1.4) is **robustly completely observable** if the identity function* $\theta(t,x) = x$ *is completely observable.*

*We say that a function* $\theta \in C^0(\Re^+ \times \Re^n ; \Re^l)$ *with* $\theta(\cdot, 0) = 0$ *is **robustly completely detectable with respect to (1.4)** if there exists a function* $\Psi \in C^0(\Re^+ \times \Re^k \times \Re^m ; \Re^l)$ *with* $\Psi(t,0,0) = 0$ *for all* $t \geq 0$, *functions* $\sigma \in KL$, $\beta \in K^+$ *such that for every* $(t_0, x_0) \in \Re^+ \times \Re^n$ *and* $d(\cdot) \in M_D$ *the solution* $x(t)$ *of (1.4) with initial condition* $x(t_0) = x_0$ *and corresponding to* $d(\cdot) \in M_D$ *satisfies:*

$$|\theta(t,x(t)) - \Psi(t, h(t,x(t)), \mathcal{D} y(t,x(t)))| \leq \sigma(\beta(t_0)|x_0|, t - t_0), \; \forall t \geq t_0 \quad (2.3)$$

*We say that system (1.4) is **robustly completely detectable** if the identity function* $\theta(t,x) = x$ *is robustly completely detectable with respect to (1.4).*

*Suppose that system (1.4) is robustly completely observable and moreover suppose that the continuous functions* $y_{(m_j+1),j} : \Re^+ \times \Re^n \to \Re$, $j = 1,...,p$, *defined below are all independent of* $d \in D$ *:*

$$y_{(m_j+1),j}(t,x) = \frac{\partial y_{m_j,j}}{\partial t}(t,x) + \frac{\partial y_{m_j,j}}{\partial x}(t,x) f(t,x,d)$$

*Define the continuous functions* $\tilde{y}_{(m_j+1),j} : \Re^+ \times \Re^k \times \Re^m \to \Re$, $j = 1,...,p$ *:*

$$\tilde{y}_{(m_j+1),j}(t,y,z) = y_{(m_j+1),j}(t, \Psi(t,y,z)) \quad (2.4)$$

*where* $\Psi \in C^0(\Re^+ \times \Re^k \times \Re^m ; \Re^n)$ *is the function for which (2.2) with* $\theta(t,x) = x$ *holds. We say that system (1.4) is **Robustly Lipschitz completely observable**, if the continuous functions* $\bar{y}_{(m_j+1),j}(t,x,z) := \tilde{y}_{(m_j+1),j}(t,h(t,x),z)$, $j = 1,...,p$ *are locally Lipschitz with respect to* $(x,z) \in \Re^n \times \Re^m$, *in the sense that, for every bounded interval* $I \subset \Re^+$ *and for every compact subset* $S \subset \Re^n \times \Re^m$, *there exists a* *constant* $L \geq 0$ *such that* $\left|\bar{y}_{(m_j+1),j}(t,x,z) - \bar{y}_{(m_j+1),j}(t,v,w)\right| = \left|\tilde{y}_{(m_j+1),j}(t,h(t,x),z) - \tilde{y}_{(m_j+1),j}(t,h(t,v),w)\right| \leq L|(x-v,z-w)|$,
*for all* $(t;(x,z);(v,w)) \in I \times S \times S$ *and* $j = 1,...,p$.

**Remark 2.5:**
**a)** If system (1.4) is robustly completely observable, then every function $\theta \in C^0(\Re^+ \times \Re^n; \Re^l)$ with $\theta(\cdot,0) = 0$ is robustly completely observable.
**b)** For a linear system $\dot{x} = A(t)x$, $y = h(t)x$, where the matrices $A(t)$, $h(t)$ have real analytic entries, complete observability is equivalent to observability (see pages 279-280 in [26]). In general, complete observability implies observability.
**c)** Notice that for every input $d(\cdot) \in M_D$ and for every $(t_0, x_0) \in \Re^+ \times \Re^n$, the unique solution $x(t)$ of (1.4) corresponding to $d(\cdot)$ and initiated from $x_0$ at time $t_0$, satisfies the following relations for $j = 1,...,p$:

$$\dot{y}_{i,j}(t) = a_j(t,y(t))y_{i+1,j}(t) - \varphi_{i+1,j}(t,y(t)), \quad \forall t \geq t_0, \; i = 0,...,m_j - 1$$

where $y_{i,j}(t) := y_{i,j}(t,x(t))$ and $y(t) = h(t,x(t))$. Thus the functions $\{\varphi_{i,j}\}$ $i = 1,...,m_j$, $j = 1,...,p$ play the role of "output injection", used in the literature for the construction of observers with linear error dynamics (see [19,25] and the references therein). If the dynamics and the output maps of system (1.4) are of class $C^\infty(\Re^+ \times \Re^n \times D; \Re^n)$ and $C^\infty(\Re^+ \times \Re^n; \Re^k)$, respectively, then the functions $\{\varphi_{i,j}\}$ $i = 1,...,m_j$, $j = 1,...,p$ can be selected to be identically equal to zero and the functions $\{g_j\}$ $j = 1,...,p$ can be selected to be equal to the identity function. **Thus the disturbance-free smooth output case, studied in the literature (see [10]) is automatically covered by Definition 2.4.** However, if the dynamics and the output maps of system (1.4) are merely locally Lipschitz and continuous, respectively, then the functions $\{\varphi_{i,j}\}$ $i = 1,...,m_j$, $j = 1,...,p$ and the functions $\{g_j\}$ $j = 1,...,p$ play a "regularizing" role. For example, consider the following single-output two-dimensional completely observable system:

$$\dot{x}_1 = |x_1| + x_2 \; ; \; \dot{x}_2 = 0 \; ; \; y = \mathrm{sgn}(x_1)|x_1|^{\frac{1}{3}}$$

Here, the output map $h(t,x) = \mathrm{sgn}(x_1)|x_1|^{\frac{1}{3}}$ is not $C^1$, but the selection $g(t,y) = y^3$ produces the smooth map $y_0(t,x) := x_1$. Moreover, $\frac{\partial y_0}{\partial t}(t,x) + \frac{\partial y_0}{\partial x}(t,x)f(t,x) = |x_1| + x_2$, which is not a $C^1$ map, but the selection $\varphi_1(t,y) = -\left|y^3\right|$ produces the smooth map $\mathcal{D} y(t,x) = y_1(t,x) = x_2$. Thus, we obtain $x = \Psi(t,h(t,x),\mathcal{D} y(t,x))$ for all $(t,x) \in \Re^+ \times \Re^2$, where $\Psi(t,y,z) := (y^3,z)$. If the functions $\{\varphi_{i,j}\}$ $i = 1,...,m_j$, $j = 1,...,p$ and the functions $\{g_j\}$ $j = 1,...,p$ were not used in Definition 2.4 then, the above system would fail to meet the requirements of complete observability.
**d)** If system (1.4) is robustly Lipschitz completely observable then the following relations are also satisfied for $j = 1,...,p$:

$$\dot{y}_{m_j,j}(t) = \tilde{y}_{(m_j+1),j}(t,y(t),\mathcal{D} y(t)), \quad \forall t \geq t_0$$

where $\mathcal{D} y(t) = \mathcal{D} y(t,x(t))$.
**e)** The problem of establishing sufficient conditions for Robust Complete Observability of a time-varying system is an open problem. However, the study of this problem is beyond the scope of the present work. In the present work our starting point is to assume Robust Complete Observability and the emphasis is placed on the design of an observer/estimator for such a system.

The following example shows that the notion of robust complete observability (detectability) allows us to consider uncertain systems with unobservable linearization.

**Example 2.6** The single-output system

$$\dot{x}_1 = x_1 + x_2^3 \; ; \; \dot{x}_2 = -x_1 x_2^2 + d(t)x_2$$
$$y = x_1, \; x = (x_1,x_2)' \in \Re^2, \; d(\cdot) \in M_{[-1,1]}$$
(2.5)

is robustly completely observable. Notice that for every $(t_0, x_0) \in \Re^+ \times \Re^2$ and $d(\cdot) \in M_{[-1,1]}$, the solution $x(t) = (x_1(t), x_2(t))$ of system (2.5) corresponding to $d(\cdot) \in M_{[-1,1]}$ with initial condition $(x_1(t_0), x_2(t_0)) = x_0$, satisfies the estimate

$$|x(t)| \leq \exp(t) |x_0|, \quad \forall t \geq t_0 \tag{2.6}$$

and consequently (by virtue of Proposition 2.3) system (2.5) is RFC. This fact follows from the evaluation of the time derivative of the function $V(x_1, x_2) = x_1^2 + x_2^2$ along the trajectories of (2.5). Specifically, we obtain $\dot{V} \leq 2V$ and inequality (2.6) is an immediate consequence. Moreover, we have:

$$\begin{pmatrix} x_1 \\ x_2 \end{pmatrix} = \Psi(t, y, \mathcal{D}y) := \begin{pmatrix} y \\ \operatorname{sgn}(\mathcal{D}y)|\mathcal{D}y|^{\frac{1}{3}} \end{pmatrix}, \; \mathcal{D}y(t,x) = y_1(t,x) := x_2^3, \; \varphi_1(t,y) = -y, \; g(t,y) := y, \; a(t,y) \equiv 1 \tag{2.7}$$

However, notice that the linearization of system (2.5) is not detectable. On the other hand the single output system

$$\dot{x}_1 = x_1 + x_2^3 \; ; \; \dot{x}_2 = -x_1 x_2^2 + d(t)x_2 \; ; \; \dot{x}_3 = -(1+|d(t)|)x_3$$
$$y = x_1 \; ; \; (x_1, x_2, x_3) \in \Re^3, \; d(\cdot) \in M_{[-1,1]} \tag{2.8}$$

is robustly completely detectable since it is RFC and the following inequality holds for the solution of (2.8) for the continuous mapping $\Psi(t, y, \mathcal{D}y) := \left( y, \; \operatorname{sgn}(\mathcal{D}y)|\mathcal{D}y|^{\frac{1}{3}}, \; 0 \right)$:

$$|x(t) - \Psi(t, y(t), \mathcal{D}y(t))| = |x_3(t)| \leq \exp(-(t-t_0))|x_3(t_0)|, \; \forall t \geq t_0$$

The notions of $\phi$–Estimator and $\phi$–Observer are crucial for the present paper. We emphasize that an estimator is not necessarily an observer since it does not necessarily satisfy the consistent initialization property (see [15]).

**Definition 2.7** Let $\phi \in K^+$ and $\theta \in C^0(\Re^+ \times \Re^n; \Re^l)$ with $\theta(\cdot, 0) = 0$. Consider system (1.4) and suppose that it is RFC. The system

$$\dot{z}(t) = k(t, z(t), y(t)) \; ; \; \overline{\theta}(t) = \Psi(t, y(t), z(t))$$
$$z \in \Re^m, t \geq 0, \overline{\theta} \in \Re^l \tag{2.9}$$

where $k \in C^0(\Re^+ \times \Re^m \times \Re^k; \Re^m)$ with $k(t, 0, 0) = 0$, the map $\tilde{k}(t, z, x) := k(t, z, h(t, x))$ is locally Lipschitz with respect to $(x, z) \in \Re^n \times \Re^m$ and $\Psi \in C^0(\Re^+ \times \Re^k \times \Re^m; \Re^l)$ with $\Psi(t, 0, 0) = 0$ for all $t \geq 0$, is called a $\phi$–**Estimator** for $\theta$ with respect to (1.4) if system (1.4) with (2.9) is RFC and there exist functions $\sigma \in KL$ and $\beta \in K^+$, such that for every $(x_0, z_0) \in \Re^n \times \Re^m$, $t_0 \geq 0$, $d(\cdot) \in M_D$, the unique solution $(x(\cdot), z(\cdot))$ of system (1.4) with (2.9) initiated from $(x_0, z_0) \in \Re^n \times \Re^m$ at time $t_0 \geq 0$ and corresponding to $d(\cdot) \in M_D$, satisfies the following estimate:

$$\phi(t)|\overline{\theta}(t) - \theta(t, x(t))| \leq \sigma\big(\beta(t_0)|(x_0, z_0)|, t - t_0\big), \; \forall t \geq t_0 \tag{2.10a}$$

System (2.9) is called a $\phi$–**Estimator** for system (1.4) if $\theta(t, x) := x$. If $\phi(t) \equiv 1$, then (2.9) is simply called an **Estimator** for $\theta$ with respect to (1.4). In any case, the continuous map $\Psi \in C^0(\Re^+ \times \Re^k \times \Re^m; \Re^l)$ is called the reconstruction map of the ($\phi$–)Estimator for $\theta$ with respect to (1.4).

**Definition 2.8** Let $\phi \in K^+$ and consider system (1.4). Suppose that (1.2) is a $\phi$-estimator for the identity function $\theta(t, x) \equiv x$ with respect to (1.4) and that (1.2) satisfies the **Consistent Initialization Property**, i.e., for every $(t_0, x_0) \in \Re^+ \times \Re^n$ there exists $z_0 \in \Re^m$ such that the solution $(x(\cdot), z(\cdot))$ of system (1.4) with (1.2) initiated from $(x_0, z_0) \in \Re^n \times \Re^m$ at time $t_0 \geq 0$ and corresponding to arbitrary $d(\cdot) \in M_D$, satisfies

$$x(t) = \Psi(t, y(t), z(t)), \; \forall t \geq t_0, \; \forall d(\cdot) \in M_D \tag{2.10b}$$

*Then we say that system (1.2) is a global $\phi$ – observer for (1.4), or that the global $\phi$ – observer problem for (1.4) is solvable. If $\phi(t) \equiv 1$ then we say that system (1.2) is a global observer for (1.4), or that the global observer problem for (1.4) is solvable. Particularly, the continuous map $\Psi \in C^0(\Re^+ \times \Re^k \times \Re^m ; \Re^n)$ is called the reconstruction map of the global observer (1.2).*

**Remark 2.9:** Necessary and sufficient conditions for the existence of a global observer for (1.4) with identity reconstruction map, i.e., $z \in \Re^n$ and $\Psi(t, y, z) \equiv z$, are given in [15], by exploiting the notion of the Observer Lyapunov Function (OLF).

**Remark 2.10:** In other words, if system (2.9) is an estimator for $\theta$ with respect to (1.4) (i.e., the case of uncertain dynamical system) then the following system

$$\dot{x}(t) = f(t, x(t), d(t)) \quad ; \quad \dot{z}(t) = k(t, z(t), h(t, x(t)))$$
$$Y(t) = \Psi(t, h(t, x(t)), z(t)) - \theta(t, x(t))$$
$$(x, z) \in \Re^n \times \Re^m , t \geq 0 , Y \in \Re^l , d(\cdot) \in M_D$$

is **non-uniformly in time Robustly Globally Asymptotically Output Stable** (RGAOS, see [14]). Moreover, there exists an estimator for $\theta$ with respect to (1.4) if and only if there exists a function $\Psi \in C^0(\Re^+ \times \Re^k \times \Re^m ; \Re^l)$ with $\Psi(t,0,0) = 0$ for all $t \geq 0$ such that the **Robust Output Feedback Stabilization problem** (ROFS problem, see [15]) with measured output $y = (h(t, x), z)$ and stabilized output $Y = \Psi(t, z, h(t, x)) - \theta(t, x)$ is globally solvable for the system:

$$\dot{x}(t) = f(t, x(t), d(t)) \quad ; \quad \dot{z}(t) = v(t)$$
$$y(t) = (h(t, x(t)), z(t)) \quad ; \quad Y(t) = \Psi(t, h(t, x(t)), z(t)) - \theta(t, x(t))$$
$$(x, z) \in \Re^n \times \Re^m , t \geq 0 , v \in \Re^m , Y \in \Re^l , d(\cdot) \in M_D , y \in \Re^k \times \Re^m$$

Consequently, by virtue of Proposition 2.6 in [15], if system (2.9) is an estimator for $\theta$ with respect to (1.4) then there exist functions $V \in C^1(\Re^+ \times \Re^n \times \Re^m ; \Re^+)$, $a_1, a_2 \in K_\infty$, $\beta, \mu \in K^+$, such that the following inequalities hold for all $(t, z, x, d) \in \Re^+ \times \Re^m \times \Re^n \times D$:

$$a_1\left(|\Psi(t, h(t, x), z) - \theta(t, x)|\right) + a_1\left(\mu(t)|(z, x)|\right) \leq V(t, z, x) \leq a_2\left(\beta(t)|(z, x)|\right)$$
$$\frac{\partial V}{\partial t}(t, z, x) + \frac{\partial V}{\partial x}(t, z, x) f(t, x, d) + \frac{\partial V}{\partial z}(t, z, x) k(t, z, h(t, x)) \leq -V(t, z, x)$$

Finally, the following technical lemmas constitute the basic tools for the proof of the existence of estimators with respect to (1.4).

**Lemma 2.11** *Let $y : \Re^+ \to \Re^+$ be an absolutely continuous function that satisfies the following differential inequality a.e. for $t \in [t_0, t_1]$:*

$$\dot{y}(t) \leq -a(t) y(t) + b(t) \tag{2.11}$$

*where $a : \Re^+ \to \Re$, $b : \Re^+ \to \Re^+$ are $C^0$ functions that satisfy $\int_0^{+\infty} a(t) dt = +\infty$, $\lim_{t \to +\infty} \frac{b(t)}{a(t)} = M \geq 0$, $a(t) > 0$ for all $t \geq T$ and for certain time $T \geq 0$. Then there exist a constant $K > 0$, being independent of $t_1$ such that:*

$$y(t) \leq K , \forall t \in [t_0, t_1] \tag{2.12}$$

**Lemma 2.12** *Consider the control system:*

$$\dot{x}_i = a(t, \theta) x_{i+1} + v_i \quad i = 1, \ldots, n-1$$
$$\dot{x}_n = v_n + u \tag{2.13}$$

*where $x = (x_1, \ldots, x_n) \in \Re^n$ is the state, $v = (v_1, \ldots, v_n) \in \Re^n$ is the control input, $\theta(t) \in \Theta \subseteq \Re^m$ is the vector of time-varying parameters and $a \in C^0(\Re^+ \times \Theta ; \Re)$ is a mapping that satisfies $\inf\{a(t, \theta) ; (t, \theta) \in \Re^+ \times \Theta\} > 0$.*

*Then for every $\phi \in K^+$ there exist $q \in K^*$, $\rho \in K^+$, a vector $k \in \Re^n$ and constants $\gamma, R, M > 0$ such that for every $(t_0, \theta, x_0, u) \in \Re^+ \times C^0(\Re^+; \Theta) \times \Re^n \times L_{loc}^\infty(\Re^+; \Re)$ the solution of the closed-loop system (2.13) with*

$$v = a(t, \theta) \, diag\big(Rq(t), R^2 q^2(t), ..., R^n q^n(t)\big) k \, x_1 \qquad (2.14)$$

*initial condition $x(t_0) = x_0$ and corresponding to inputs $(\theta, u) \in C^0(\Re^+; \Theta) \times L_{loc}^\infty(\Re^+; \Re)$ satisfies the estimate for all $t \geq t_0$:*

$$\phi(t)|x(t)| \leq \rho(t_0) \exp\big(-\gamma(t - t_0)\big)|x_0| + M \sup_{\tau \in [t_0, t]} \left( \frac{|u(\tau)|}{\phi(\tau)} \right) \qquad (2.15)$$

**Lemma 2.13** *Suppose that $\Psi \in C^0(\Re^+ \times \Re^m; \Re^n)$. Then there exist functions $a_1, a_2 \in K_\infty$, $\beta \in K^+$ such that for every $(t, x, y) \in \Re^+ \times \Re^m \times \Re^m$, it holds that:*

$$|\Psi(t, x) - \Psi(t, y)| \leq a_1\big(\beta(t)|x - y|\big) + a_1\big(a_2(|y|)|x - y|\big) \qquad (2.16)$$

## 3. Main Results and Examples

Our first main result states that we may design an estimator for a robustly completely detectable function. Moreover, we can assign the convergence rate under the hypothesis of robust complete observability.

**Theorem 3.1** *If the function $\theta \in C^0(\Re^+ \times \Re^n; \Re^l)$ with $\theta(\cdot, 0) = 0$ is robustly completely observable with respect to (1.4) then for every $\phi \in K^+$ there exists a $\phi$-estimator for $\theta$ with respect to (1.4). If the function $\theta \in C^0(\Re^+ \times \Re^n; \Re^l)$ is robustly completely detectable then there exists an estimator for $\theta$ with respect to (1.4).*

**Proof** For simplicity of the proof we assume that $p = 1$. The case $p > 1$ is completely analogous. Let $\phi \in K^+$ be given. Clearly, the hypotheses made guarantee the existence of a function $\Psi \in C^0(\Re^+ \times \Re^k \times \Re^m; \Re^l)$ with $\Psi(t, 0, 0) = 0$ for all $t \geq 0$ such that

$$\theta(t, x) = \Psi(t, h(t, x), \mathcal{D} y(t, x)), \quad \forall (t, x) \in \Re^+ \times \Re^n \qquad (3.1)$$

(case of robust complete observability) or there exist functions $\sigma \in KL$, $\beta \in K^+$ such that for every $(t_0, x_0) \in \Re^+ \times \Re^n$ and $d(\cdot) \in M_D$ the solution $x(t)$ of (1.4) with initial condition $x(t_0) = x_0$ and corresponding to $d(\cdot) \in M_D$ satisfies:

$$\big|\theta(t, x(t)) - \Psi(t, h(t, x(t)), \mathcal{D} y(t, x(t)))\big| \leq \sigma\big(\beta(t_0)|x_0|, t - t_0\big), \quad \forall t \geq t_0 \qquad (3.2)$$

(case of robust complete detectability), where

$$\mathcal{D} y(t, x) := (y_1(t, x), ..., y_m(t, x)) \in \Re^m \qquad (3.3)$$

and the functions $y_i(t, x)$ are defined recursively by the following relations for $i = 0, ..., m$:

$$y_0(t, x) = g(t, h(t, x))$$

$$y_i(t, x) := \frac{1}{a(t, h(t, x))} \left\{ \frac{\partial y_{i-1}}{\partial t}(t, x) + \frac{\partial y_{i-1}}{\partial x}(t, x) f(t, x, d) + \varphi_i(t, h(t, x)) \right\}, \quad i = 1, ..., m$$

and are all independent of $d \in D$ and of class $C^1(\Re^+ \times \Re^n; \Re)$, where $g : \Re^+ \times \Re^k \to \Re$, $a : \Re^+ \times \Re^k \to \Re$ are functions of class $C^0(\Re^+ \times \Re^k; \Re)$ and $\{\varphi_i\}$ $i = 1, ..., m$ are functions of class $C^0(\Re^+ \times \Re^k; \Re)$ with $\inf\{a(t, y); (t, y) \in \Re^+ \times \Re^k\} > 0$, $g(t, 0) = \varphi_i(t, 0) = 0$ for all $t \geq 0$. Clearly, the above definitions guarantee that $\mathcal{D} y(t,0) = 0$ for all $t \geq 0$ and using Fact V in [13] in conjunction with the regularity properties of $f, g, a, h, \mathcal{D} y$, we obtain functions $p \in K_\infty$ and $\beta \in K^+$ such that:

$$|h(t,x)| + |g(t,h(t,x))| + |\mathcal{D} y(t,x)| + \left| \frac{\partial y_m}{\partial t}(t,x) + \frac{\partial y_m}{\partial x}(t,x) f(t,x,d) \right| \leq p(\beta(t)|x|), \ \forall (t,x,d) \in \Re^+ \times \Re^n \times D \quad (3.4)$$

Since, system (1.4) is RFC, by virtue of Proposition 2.3, there exist functions $\tilde{a} \in K_\infty$ and $\mu \in K^+$ such that for every input $d(\cdot) \in M_D$ and for every $(t_0, x_0) \in \Re^+ \times \Re^n$, the unique solution $x(t)$ of (1.4) corresponding to $d(\cdot)$ and initiated from $x_0$ at time $t_0$ exists for all $t \geq t_0$ and satisfies:

$$|x(t)| \leq \mu(t) \tilde{a}(|x_0|), \ \forall t \geq t_0 \quad (3.5)$$

Moreover, by virtue of Lemma 2.13, there exist functions $a_1, a_2 \in K_\infty$, $\gamma \in K^+$ such that for every $(t, y, z, \mathcal{D} y) \in \Re^+ \times \Re^k \times \Re^m \times \Re^m$, it holds that:

$$|\Psi(t, y, z) - \Psi(t, y, \mathcal{D} y)| \leq a_1(\gamma(t)|z - \mathcal{D} y|) + a_1(a_2(|y| + |\mathcal{D} y|)|z - \mathcal{D} y|) \quad (3.6)$$

Notice that for every input $d(\cdot) \in M_D$ and for every $(t_0, x_0) \in \Re^+ \times \Re^n$, the components of the vector $(y_1(t), ..., y_m(t)) = \mathcal{D} y(t) = \mathcal{D} y(t, x(t)) := (y_1(t, x(t)), ..., y_m(t, x(t))) \in \Re^m$, where $x(t)$ denotes the unique solution of (1.4) corresponding to $d(\cdot)$ and initiated from $x_0$ at time $t_0$, satisfy the following relations:

$$\dot{y}_i(t) = a(t, y(t)) y_{i+1}(t) - \varphi_{i+1}(t, y(t)), \ \forall t \geq t_0, \ i = 0, ..., m-1 \quad (3.7)$$

where $y(t) = h(t, x(t))$ and $y_0(t) := g(t, y(t))$. By virtue of Corollary 10 and Remark 11 in [27], there exists a function $\kappa \in K_\infty$ such that

$$p(r s) + a_1(r s) + a_2(r s) \leq \kappa(r) \kappa(s), \ \forall r, s \geq 0 \quad (3.8)$$

where $p \in K_\infty$ is the function involved in (3.4) and $a_1, a_2 \in K_\infty$ are the functions involved in (3.6). It follows by (3.4), (3.5), (3.7) and (3.8) that for every input $d(\cdot) \in M_D$ and for every $(t_0, x_0) \in \Re^+ \times \Re^n$, we have:

$$|y(t)| + |y_0(t)| + |\mathcal{D} y(t)| + |\dot{y}_m(t)| \leq \tilde{\beta}(t) \kappa(\tilde{a}(|x_0|)), \ \forall t \geq t_0 \quad (3.9)$$

where $\tilde{\beta} \in K^+$ is a function that satisfies $\tilde{\beta}(t) \geq \kappa(\beta(t)\mu(t))$ for all $t \geq 0$. Let $\tilde{\phi} \in K^+$ a function that satisfies:

$$\tilde{\phi}(t) \geq \tilde{\beta}(t) + \frac{\gamma(t) + \kappa(\tilde{\beta}(t))}{\kappa^{-1}\left(\frac{\exp(-t)}{\phi(t)}\right)}, \ \forall t \geq 0 \quad (3.10)$$

where $\phi \in K^+$ is the given function of our problem and $\gamma \in K^+$ is the function involved in (3.6). By virtue of Lemma 2.12 there exist $q \in K^*$, $\rho \in K^+$, a vector $k \in \Re^{m+1}$ and constants $R, M > 0$ such that for every $(t_0, z(t_0), \dot{y}_m) \in \Re^+ \times \Re^{m+1} \times L^\infty_{loc}(\Re^+; \Re)$ the solution of the closed-loop system

$$\begin{aligned} \dot{z}_i &= a(t, y(t)) z_{i+1} - \varphi_i(t, y(t)) + v_i \quad i = 0, ..., m-1 \\ \dot{z}_m &= v_m \\ \bar{\theta} &:= \Psi(t, y, z_1, ..., z_m) \ ; \ z := (z_0, ..., z_m) \in \Re^{m+1}, \ v := (v_0, ..., v_m) \in \Re^{m+1} \end{aligned} \quad (3.11)$$

with

$$v = a(t, y(t)) \, diag(Rq(t), R^2 q^2(t), ..., R^{m+1} q^{m+1}(t)) k (z_0 - g(t, y(t))) \quad (3.12)$$

initial condition $z(t_0)$ and corresponding to input $\dot{y}_m \in L^\infty_{loc}(\Re^+; \Re)$ satisfies the estimate for all $t \geq t_0$:

$$\tilde{\phi}(t) |z(t) - (g(t, y(t)), \mathcal{D} y(t))| \leq \rho(t_0) \exp(-\gamma(t - t_0)) |z(t_0) - (g(t_0, y(t_0)), \mathcal{D} y(t_0))| + M \sup_{\tau \in [t_0, t]} \left( \frac{|\dot{y}_m(\tau)|}{\tilde{\phi}(\tau)} \right) \quad (3.13)$$

Clearly, inequalities (3.4), (3.9), (3.10) in conjunction with (3.13) imply that:

$$\tilde{\phi}(t)|z(t) - (g(t, y(t)), \mathcal{D} y(t))| \leq \rho(t_0)|z(t_0)| + \rho(t_0)p(\beta(t_0)|x(t_0)|) + M \kappa(\tilde{a}(|x(t_0)|)), \quad \forall t \geq t_0 \quad (3.14)$$

Making use of inequalities (3.6), (3.8), (3.9), (3.10) and (3.14) we conclude:

$$\phi(t)|\Psi(t, y(t), z_1(t),..., z_m(t)) - \Psi(t, y(t), \mathcal{D} y(t))| \leq \exp(-(t-t_0))\omega(\lambda(t_0)|(x(t_0), z(t_0))|), \quad \forall t \geq t_0 \quad (3.15)$$

for appropriate functions $\omega \in K_\infty$ and $\lambda \in K^+$. Notice that, by virtue of (3.5), (3.9) and (3.14), system (1.4) with (3.11) is RFC. Inequality (3.15) in conjunction with (3.1) or (3.2), proves the statement of the theorem. ◁

**Example 3.2** Consider again system (2.5), where the mappings $\mathcal{D} y(t, x) = y_1(t, x) := x_2^3$, $\varphi_1(t, y) = -y$, $g(t, y) := y$ and $(x_1, x_2) = \Psi(t, y, \mathcal{D} y)$ are given in (2.7). Notice that by virtue of (2.6) we have:

$$|y(t)| + |\mathcal{D} y(t)| + |\dot{y}_1(t, x(t))| \leq 3 \exp(5t) \left(|x_0| + |x_0|^3 + |x_0|^5\right), \quad \forall t \geq t_0 \quad (3.16)$$

Making use of the inequality $\left|\text{sgn}(x)|x|^{\frac{1}{3}} - \text{sgn}(y)|y|^{\frac{1}{3}}\right| \leq 2|x-y|^{\frac{1}{3}}$, which holds for all $x, y \in \Re$, we obtain that:

$$|\Psi(t, y, z) - \Psi(t, y, \mathcal{D} y)| \leq 2|z - \mathcal{D} y|^{\frac{1}{3}}, \quad \forall (t; y, z, \mathcal{D} y) \in \Re^+ \times \Re^3 \quad (3.17)$$

Using Lemma 2.12 for $\phi(t) = \exp(5t)$ we obtain that there exists a function $\rho \in K^+$ and constants $M, R > 0$

$$\begin{aligned}
\dot{z}_1 &= y(t) + z_2 - 12R \exp(10t)(z_1 - y(t)) \quad ; \quad \dot{z}_2 = -72R^2 \exp(20t)(z_1 - y(t)) \\
\bar{x} &= \Psi(t, y, z_2) \quad ; \quad z := (z_1, z_2) \in \Re^2, \, t \geq 0
\end{aligned} \quad (3.18)$$

such that for every $(t_0, z(t_0)) \in \Re^+ \times \Re^2$, the solution of (3.18) satisfies:

$$|z(t) - (y(t), \mathcal{D} y(t))| \leq \exp(-5t) \left(\rho(t_0)|z(t_0) - (y(t_0), \mathcal{D} y(t_0))| + M \left(|x(t_0)| + |x(t_0)|^3 + |x(t_0)|^5\right)\right), \quad \forall t \geq t_0 \quad (3.19)$$

Let $a(s) := 2s^{\frac{1}{3}}$. It follows from (3.17) that we have:

$$|x(t) - \Psi(t, y(t), z_2(t))| \leq \exp(-(t-t_0))a\left(\rho(t_0)|z(t_0) - (y(t_0), \mathcal{D} y(t_0))| + M\left(|x(t_0)| + |x(t_0)|^3 + |x(t_0)|^5\right)\right), \quad \forall t \geq t_0 \quad (3.20)$$

Thus we may conclude that system (3.18) is an estimator for system (2.5), which guarantees exponential convergence. ◁

The following theorem deals with the solvability of the global $\phi$ – observer problem for (1.4).

**Theorem 3.3** *Suppose that (1.4) is Robustly Lipschitz completely observable with $p = 1$. Then the global $\phi$ – observer problem for (1.4) is solvable for all $\phi \in K^+$.*

**Proof** Let $\phi \in K^+$ be arbitrary. The same arguments in the proof of Theorem 4.1 are repeated (mathematical relations (3.1)-(3.9)) with $\theta(t, x) \equiv x$ and we define

$$y_{m+1}(t, x) := \frac{\partial y_m}{\partial t}(t, x) + \frac{\partial y_m}{\partial x}(t, x) f(t, x, d) \quad (3.21a)$$

then it is a consequence of (2.4) that the following equality holds for all $(t, x) \in \Re^+ \times \Re^n$:

$$\tilde{y}_{m+1}(t, h(t, x), \mathcal{D} y(t, x)) = y_{m+1}(t, x) \quad (3.21b)$$

Notice that by (3.21a,b) and (3.1), (3.7) we obtain: $\tilde{y}_{m+1}(t, y(t), \mathcal{D} y(t)) = \dot{y}_m(t)$, for all $t \geq t_0$ and $d \in M_D$. Moreover, since (1.4) is Robustly Lipschitz completely observable, it follows that the function $\bar{y}_{m+1}(t, x, z) := \tilde{y}_{m+1}(t, h(t, x), z)$ is locally Lipschitz with respect to $(x, z) \in \Re^n \times \Re^m$. Define the function

$$\beta(t,w) := \tilde{\beta}(t)(1+\exp(t)|w|) \qquad (3.22)$$

where $\tilde{\beta} \in K^+$ is the function involved in (3.9). Notice that by virtue of (3.9) and definition (3.21), we have

$$|\tilde{y}_m(t)| \le \beta(t,w(t)), \quad \forall t \ge t_0, \, d \in M_D, \text{ provided that } \kappa(\tilde{a}(|x_0|)) \le 1+\exp(t)|w(t)|, \text{ for all } t \ge t_0 \qquad (3.23)$$

Let $\tilde{\phi} \in K^+$ be a function that satisfies (3.10), where $\gamma \in K^+$ is the function involved in (3.6) and $\phi \in K^+$ is the given function of our problem. By virtue of Lemma 2.12 there exist $q \in K^*$, $\rho \in K^+$, a vector $k \in \Re^{m+1}$ and constants $R, M > 0$ such that for every $(t_0, z(t_0), w(t_0), \tilde{y}_m) \in \Re^+ \times \Re^{m+1} \times \Re \times L_{loc}^\infty(\Re^+; \Re)$ the solution of the closed-loop system

$$\begin{aligned}
\dot{z}_i &= a(t,y(t))z_{i+1} - \varphi_i(t,y(t)) + v_i \quad i = 0,\ldots,m-1 \\
\dot{z}_m &= \beta(t,w)\,sat\!\left(\frac{\tilde{y}_{m+1}(t,y(t),z_1,\ldots,z_m)}{\beta(t,w)}\right) + v_m \\
\dot{w} &= -w \\
\bar{x} &:= \Psi(t,y,z_1,\ldots,z_m), \, z := (z_0,\ldots,z_m) \in \Re^{m+1}, \, v := (v_0,\ldots,v_m) \in \Re^{m+1}
\end{aligned} \qquad (3.24)$$

with (3.12) satisfies the following estimate for all $t \ge t_0$, $d \in M_D$ (recall that $y_0(t) = g(t,y(t))$):

$$\tilde{\phi}(t)|z(t) - (y_0(t), \mathcal{D}y(t))|$$
$$\le \rho(t_0)|z(t_0) - (g(t_0,y(t_0)), \mathcal{D}y(t_0))| + M \sup_{t_0 \le \tau \le t} \left| \frac{\tilde{y}_m(\tau)}{\tilde{\phi}(\tau)} - \frac{\beta(\tau,w(\tau))}{\tilde{\phi}(\tau)} sat\!\left(\frac{\tilde{y}_{m+1}(\tau,y(\tau),z_1(\tau),\ldots,z_m(\tau))}{\beta(\tau,w(\tau))}\right) \right| \qquad (3.25)$$

First, notice that, if $z(t_0) = (g(t_0,y(t_0)), \mathcal{D}y(t_0))$ and $\kappa(\tilde{a}(|x(t_0)|)) \le \exp(t_0)|w(t_0)|$, then it follows from (3.23) (and the fact that $w(t) = \exp(-(t-t_0))w(t_0)$ for all $t \ge t_0$) that $z(t) = (y_0(t), \mathcal{D}y(t))$ for all $t \ge t_0$, $d \in M_D$ and consequently, by virtue of (3.1), $\bar{x}(t) = x(t)$ for all $t \ge t_0$, $d \in M_D$. Moreover, notice that by virtue of definition (3.22), inequalities (3.4), (3.9), (3.10) and (3.25) we obtain the following inequality for all $t \ge t_0$ and $d \in M_D$:

$$\tilde{\phi}(t)|z(t) - (y_0(t), \mathcal{D}y(t))| \le \rho(t_0)|z(t_0)| + \rho(t_0)p(\beta(t_0)|x(t_0)|) + M\,\kappa(\tilde{a}(|x(t_0)|)) + M(1+|w(t_0)|) \qquad (3.26)$$

which in conjunction with estimates (3.5) and (3.9) (and the fact that $w(t) = \exp(-(t-t_0))w(t_0)$ for all $t \ge t_0$), shows that (1.4) with (3.24) is RFC (recall Proposition 2.3). Making use of inequalities (3.6), (3.8), (3.9), (3.10) and (3.13) we conclude:

$$\phi(t)|\Psi(t,y(t),z_1(t),\ldots,z_m(t)) - \Psi(t,y(t),\mathcal{D}y(t))| \le \exp(-(t-t_0))\left(\omega(\lambda(t_0)|(x(t_0),z(t_0),w(t_0))|) + Q\right), \, \forall t \ge t_0, \, d \in M_D \qquad (3.27)$$

for appropriate constant $Q > 0$ and functions $\omega \in K_\infty$ and $\lambda \in K^+$. It follows from Lemma 3.5 in [14], that since:

(i) Estimate (3.27) and equation (3.1) with $\theta(t,x) \equiv x$ hold (which guarantees Robust Global Output Attractivity for the output $Y := \phi(t)(\Psi(t,h(t,x),z_1,\ldots,z_m) - x)$,

(ii) System (1.4) with (3.24) is RFC

(iii) The point $(x,z,w) = (0,0,0)$ is the equilibrium point of (1.4) with (3.24),

(iv) The dynamics of (1.4) with (3.24) are locally Lipschitz with respect to $(x,z,w) \in \Re^n \times \Re^{m+1} \times \Re$

that system (1.4) with (3.24) is non-uniformly in time Robustly Globally Asymptotically Output Stable (RGAOS, see [14]), with output defined by $Y := \phi(t)(\Psi(t,h(t,x),\xi_2,\ldots,\xi_{m+1}) - x)$. Consequently, by virtue of Lemma 3.4 in [14], it follows that there exist functions $\sigma \in KL$ and $b \in K^+$ such that for every $t_0 \ge 0$ and $(x(t_0), z(t_0), w(t_0)) \in \Re^n \times \Re^{m+1} \times \Re$ the following estimate holds:

$$\phi(t)|\Psi(t,y(t),z_1(t),\ldots,z_m(t)) - x(t)| \le \sigma(b(t_0)|(x(t_0),z(t_0),w(t_0))|, t-t_0), \, \forall t \ge t_0 \qquad (3.28)$$

Estimate (3.28) implies that (3.24) is a $\phi$-estimator for the identity function $\theta(t,x) \equiv x$ with respect to (1.4). Since (3.24) is a $\phi$-estimator for system (1.4), which satisfies the consistent initialization property and since $\phi \in K^+$ is arbitrary, we conclude that the global $\phi$-observer problem for (1.4) is solvable. ◁

Notice that the constructed observer (3.24) in the above proof of Theorem 3.3 is a high-gain type observer (with increasing time-varying gains). High-gain type observers were also considered in [8-10] for autonomous systems with analytic dynamics. The following example illustrates that Theorem 3.3 is applicable even to **autonomous disturbance-free systems with smooth outputs and unobservable linearization**.

**Example 3.4** Consider again system (2.5), where $d \in C^0(\Re^+;[-1,1])$ is a known function. Particularly, when $d \in C^0(\Re^+;[-1,1])$ is a constant function then system (2.5) is an autonomous system with unobservable linearization. In this case we have:

$$y = y_0(t,x) = x_1 \; ; \; \mathcal{D}\, y(t,x) = y_1(t,x) = x_2^3$$

$$\varphi_1(t,y) := -y \, , \; g(t,y) := y \, , \; a(t,y) \equiv 1 \; ; \; y_2(t,x) = 3d(t)x_2^3 - 3x_1 x_2^4$$

$$x = \Psi(t,y,\mathcal{D}\,y) := \begin{pmatrix} y \\ \operatorname{sgn}(\mathcal{D}\,y)|\mathcal{D}\,y|^{\frac{1}{3}} \end{pmatrix} \; ; \; \tilde{y}_2(t,y,z) = 3d(t)z - 3y|z|^{\frac{4}{3}}$$

It is clear that system (2.5) with $d \in C^0(\Re^+;[-1,1])$ being a known function is Lipschitz completely observable. Since (3.16) still holds, Theorem 3.3 guarantees that there exists $R > 0$ such that the system

$$\dot{z}_1(t) = y(t) + z_2(t) - 12R \exp(10t)(z_1(t) - y(t))$$

$$\dot{z}_2(t) = \beta(t)(1 + \exp(t)|z_3(t)|) \operatorname{sat}\left( \frac{3d(t)z_2(t) - 3y(t)|z_2(t)|^{\frac{4}{3}}}{\beta(t)(1+\exp(t)|z_3(t)|)} \right) - 72R^2 \exp(20t)(z_1(t) - y(t))$$

$$\dot{z}_3(t) = -z_3(t)$$

$$\bar{x} = \Psi(t,y,z_2) \, , \; z := (z_1, z_2, z_3) \in \Re^3 \, , \; t \geq 0$$

where $\beta(t) := 3\exp(5t)$, is a global observer for system (2.5), which guarantees exponential convergence. ◁

## 5. Conclusions

In this paper we have given sufficient conditions for the existence of estimators and the solvability of the global observer problem for dynamical systems. It is showed that if a time-varying uncertain system is robustly completely detectable then there exists an estimator for this system, i.e. we can estimate asymptotically the state vector of the system. Moreover, if a time-varying uncertain system is robustly completely observable then there exists an estimator for this system that guarantees convergence of the estimates with "arbitrary fast" rate of convergence. Finally, it is proved that under the assumption of Robust Lipschitz complete observability, there is a global solution of the observer problem for a time-varying system.

# Appendix

**Proof of Lemma 2.11:** Clearly, the differential inequality (2.11) implies that:

$$y(t) \leq \exp\left(-\int_{t_0}^{t} a(s)ds\right) y(t_0) + \int_{t_0}^{t} \exp\left(-\int_{\tau}^{t} a(s)ds\right) b(\tau)d\tau, \quad \forall t \in [t_0, t_1] \tag{A1}$$

Moreover, since $a(t) > 0$ for all $t \geq T$, we obtain for all $t_0 \geq 0$ and $t \geq t_0$:

$$\int_{t_0}^{t} a(s)ds = \int_{t_0}^{t} |a(s)|ds + 2\int_{t_0}^{t} \min\{0, a(s)\}ds \geq 2\int_{t_0}^{t} \min\{0, a(s)\}ds \geq 2\int_{0}^{t} \min\{0, a(s)\}ds \geq 2\int_{0}^{T} \min\{0, a(s)\}ds$$

We define $K_1 := \exp\left(-2\int_0^T \min\{0, a(s)\}ds\right) y(t_0)$ and the previous inequalities in conjunction with (A1) give:

$$y(t) \leq K_1 + \exp\left(-\int_0^t a(s)ds\right)\int_0^t \exp\left(\int_0^\tau a(s)ds\right) b(\tau)d\tau, \quad \forall t \in [t_0, t_1] \tag{A2}$$

We define the function $p(t) := \int_0^t \exp\left(\int_0^\tau a(s)ds\right) b(\tau)d\tau$. This function is non-decreasing and consequently we either have $p(t) \leq K_3$ for some $K_3 > 0$ or $\lim_{t \to +\infty} p(t) = +\infty$. For the first case inequality (2.12) is implied by (A2) with $K = K_1 + K_1 K_3$. For the second case notice that since $\int_0^{+\infty} a(t)dt = +\infty$ and $\lim_{t \to +\infty} \frac{b(t)}{a(t)} = M \geq 0$, we can apply L'Hospital's rule for the function $q(t) := p(t)\exp\left(-\int_0^t a(s)ds\right)$ and obtain that $\lim_{t \to +\infty} q(t) = M$. Thus we may define $K := K_1 + \sup_{t \geq 0} q(t)$ and (2.12) is implied by inequality (A2). The proof is complete. ◁

**Proof of Lemma 2.12:** Without loss of generality we may assume that the given function of our problem $\phi \in K^+$ is continuously differentiable (if $\phi \in K^+$ is not continuously differentiable we may replace it by a function $\tilde{\phi} \in K^+$ which satisfies $\tilde{\phi}(t) \geq \phi(t)$ for all $t \geq 0$). Let $q \in K^*$ a function that satisfies:

$$q(t) \geq |\dot{\phi}(t)|\phi^{-1}(t) + \phi^2(t), \quad \forall t \geq 0 \tag{A3}$$

where $\phi \in K^+$ is the given function of our problem. Let $A := \{a_{i,j}; i, j = 1,...,n\}$ with $a_{i,j} := 1$ if $j = i+1$, $i = 1,...,n-1$, $a_{i,j} = 0$ if otherwise and let $c' = (1,0,...,0) \in \Re^n$. There exist a vector $k = (k_1,...,k_n)' \in \Re^n$, constants $\mu, K_1, K_2 > 0$ and a positive definite symmetric matrix $P \in \Re^{n \times n}$, such that:

$$P(A + kc') + (A + kc')'P \leq -\mu P \tag{A4}$$

$$K_1 I \leq P \leq K_2 I \tag{A5}$$

where $I \in \Re^{n \times n}$ denotes the identity matrix. Let

$$l := \inf\{a(t, \theta); (t, \theta) \in \Re^+ \times \Theta\} > 0 \tag{A6}$$

Define:

$$R := \max\left\{1; \frac{8\sqrt{n}K_2}{\mu K_1 l}\right\} \tag{A7}$$

Consider the time-varying transformation:

$$R^{i-1} q^{i-1}(t) y_i = \phi(t) x_i, \quad i = 1,...,n \tag{A8}$$

Define:

$$\tilde{F}_i(t, y_i) := \frac{\dot{\phi}(t)}{R\phi(t)q(t)} y_i, \quad i = 1,...,n \tag{A9}$$

It follows from (A2) and definition (A9) that the following inequalities hold for all $t \geq 0$ and $y_i \in \Re$:

$$|\tilde{F}_i(t, y_i)| \leq \frac{\mu l}{8\sqrt{n}} \frac{K_1}{K_2} |y_i|, \quad i = 2,...,n \tag{A10}$$

For every input $(\theta, u) \in C^0(\Re^+; \Theta) \times L^\infty_{loc}(\Re^+; \Re)$ the solution of the closed-loop system (2.13) with (2.14) is described in $y$-coordinates by the following system of differential equations:

$$\dot{y} = a(t, \theta)Rq(t)(A + kc')y + Rq(t)\tilde{F}(t, y) - \dot{q}(t)q^{-1}(t)B\,y + \phi(t)R^{1-n}q^{1-n}(t)bu \tag{A11}$$

where $y = (y_1, ..., y_n) \in \Re^n$, $B := diag(0,1,...,n-1)$, $\tilde{F}(t, y) := \left(\tilde{F}_1(t, y_1), \tilde{F}_2(t, y_2), ..., \tilde{F}_n(t, y_n)\right)'$ (defined by (A9)) and $b := (0,...,0,1)'$.

Let arbitrary $(t, \theta, x_0, u) \in \Re^+ \times C^0(\Re^+; \Theta) \times \Re^n \times L^\infty_{loc}(\Re^+; \Re)$ and consider the solution $x(t)$ of the closed-loop system (2.13) with (2.14), initial condition $x(t_0) = x_0$ and corresponding to inputs $(\theta, u) \in C^0(\Re^+; \Theta) \times L^\infty_{loc}(\Re^+; \Re)$. Clearly, for the solution $x(t)$ there exists a maximal existence time $t_{\max} > t_0$ such that the solution is defined on $[t_0, t_{\max})$ and cannot be further continued. Define the function $V(t) = y'(t)Py(t)$, where $y(t)$ is defined by the transformation (A8). By virtue of (A4), (A5), (A6), (A10) and (A11) the derivative of $V(t)$ satisfies for all $t \in [t_0, t_{\max})$ except of a set of zero Lebesgue measure:

$$\dot{V}(t) \leq -\mu l Rq(t)V(t) + 2Rq(t)K_2|y(t)|\|\tilde{F}(t, y(t))\| + 2K_2|y(t)|^2 \frac{\dot{q}(t)}{q(t)}|B| + 2K_2|y(t)|\frac{\phi(t)}{R^{n-1}q^{n-1}(t)}|u(t)|$$

$$\leq -Rq(t)\left(\frac{\mu}{4}l - \frac{2K_2}{K_1R}(n-1)\frac{\dot{q}(t)}{q^2(t)}\right)V(t) + \frac{4K_2^2}{\mu K_1 l}Rq(t)\frac{\phi^2(t)|u(t)|^2}{R^{2n}q^{2n}(t)}$$

The above differential inequality in conjunction with (A5) and inequality $|y| \leq \phi(t)|x| \leq R^{n-1}q^{n-1}(t)|y|$, which is a direct implication of definition (A8), implies that the following estimate for all $t \in [t_0, t_{\max})$:

$$\phi^2(t)|x(t)|^2 \leq \frac{K_2}{K_1} R^{2(n-1)} q^{2(n-1)}(t) \left(\frac{q(t)}{q(t_0)}\right)^{\frac{2K_2(n-1)}{K_1}} \exp\left(-\frac{\mu R l}{4}\int_{t_0}^t q(s)ds\right)\phi^2(t_0)|x_0|^2$$

$$+ \frac{4K_2^2}{\mu K_1^2 l} \int_{t_0}^t \left(\frac{q(t)}{q(\tau)}\right)^{2(n-1)\left(1+\frac{K_2}{K_1}\right)} \exp\left(-\frac{\mu R l}{4}\int_\tau^t q(s)ds\right)\frac{\phi^2(\tau)|u(\tau)|^2}{Rq(\tau)}d\tau$$

The above estimate implies that the solution $x(t)$ of the closed-loop system (2.13) with (2.14), initial condition $x(t_0) = x_0$ and corresponding to inputs $(\theta, u) \in C^0(\Re^+; \Theta) \times L^\infty_{loc}(\Re^+; \Re)$ exists for all $t \geq t_0$ (i.e., $t_{\max} = +\infty$). Define:

$$I(t) := \frac{16K_2^2}{\mu^2 K_1^2 R^2 l^2} \int_0^t \left(\frac{q(t)}{q(\tau)}\right)^{2(n-1)\left(1+\frac{K_2}{K_1}\right)} \exp\left(-\frac{\mu R l}{4}\int_\tau^t q(s)ds\right)\frac{\mu R l}{4}q(\tau)d\tau \tag{A12}$$

Clearly, the estimate given above for the solution $x(t)$ in conjunction with definition (A12) implies the following estimate for all $t \geq t_0$:

$$\phi(t)|x(t)| \leq \left(\frac{K_2}{K_1}\right)^{\frac{1}{2}} R^{n-1} q^{n-1}(t_0) \left(\frac{q(t)}{q(t_0)}\right)^{(n-1)\left(1+\frac{K_2}{K_1}\right)} \exp\left(-\frac{\mu r}{8}l\int_{t_0}^t q(s)ds\right)\phi(t_0)|x_0| + (I(t))^{\frac{1}{2}} \sup_{\tau \in [t_0, t]}\left(\frac{\phi(\tau)|u(\tau)|}{q(\tau)}\right)$$

(A13)

Integrating by parts and using inequalities (2a), we obtain:

$$I(t) \leq \frac{16K_2^2}{\mu^2 K_1^2 R^2 l^2}\left(1 + 2(n-1)\left(1+\frac{K_2}{K_1}\right)g(t)\right) \tag{A14}$$

where

$$g(t) := \int_0^t \left(\frac{q(t)}{q(\tau)}\right)^{2(n-1)\left(1+\frac{K_2}{K_1}\right)} \exp\left(-\frac{\mu R}{4}l\int_\tau^t q(s)ds\right)\frac{\dot{q}(\tau)}{q(\tau)}d\tau \tag{A15}$$

Definition (A15) implies that $g(t)$ satisfies the following differential equation:

$$\dot{g}(t) = -q(t)\left(\frac{\mu R}{4}l - 2(n-1)\left(1+\frac{K_2}{K_1}\right)\dot{q}(t)q^{-2}(t)\right)g(t) + \dot{q}(t)q^{-1}(t) \tag{A16}$$

Lemma 2.11 in conjunction with the fact that $q \in K^*$ (and thus $\lim_{t \to +\infty} \dot{q}(t)q^{-2}(t) = 0$) implies that there exists a constant $G > 0$ such that $g(t) \leq G$ for all $t \geq 0$. It follows from estimate (A13) and inequalities (A3), (A14) that there exist a constant $M_2 > 0$ such that:

$$\phi(t)|x(t)| \leq M_1 q^{n-1}(t_0) \left(\frac{q(t)}{q(t_0)}\right)^a \exp\left(-2\gamma \int_{t_0}^{t} q(s)ds\right) \phi(t_0)|x_0| + M_2 \sup_{\tau \in [t_0,t]} \left(\frac{|u(\tau)|}{\phi(\tau)}\right) \tag{A17}$$

for all $t \geq t_0$ with $M_1 := \left(\frac{K_2}{K_1}\right)^{\frac{1}{2}} R^{n-1}$, $a := (n-1)\left(1 + \frac{K_2}{K_1}\right)$ and $\gamma := \frac{\mu R}{16} l$. Next consider the function

$$y(t) := q^a(t) \exp\left(-\gamma \int_0^t q(s)ds\right) \tag{A18}$$

It is clear that $y(t)$ satisfies the differential equation $\dot{y}(t) = -q(t)\left(\gamma - a\frac{\dot{q}(t)}{q^2(t)}\right) y(t)$. Since $q \in K^*$ (which implies that $\lim_{t \to +\infty} \dot{q}(t)q^{-2}(t) = 0$), by virtue of Lemma 2.11, there exists a constant $K > 0$ such that $y(t) \leq K$ for all $t \geq 0$. Combining estimate (A17) with definition (A18) and using the fact that $q(t) \geq 1$ for all $t \geq 0$, we obtain:

$$\phi(t)|x(t)| \leq M_1 K \phi(t_0) \exp(-\gamma(t-t_0))|x_0| + M_2 \sup_{\tau \in [t_0,t]} \left(\frac{|u(\tau)|}{\phi(\tau)}\right) \tag{A19}$$

It is clear from (A19) that estimate (2.15) holds with $M := M_2$ and $\rho(t) := M_1 K \phi(t) \exp\left(\gamma \int_0^t q(s)ds\right)$. ◁

**Proof of Lemma 2.13** Clearly the function

$$\gamma(r,s) := \sup\left\{ |\Psi(t,x) - \Psi(t,y)| \, ; \, |(t,y)| \leq r, |x-y| \leq s \right\} \tag{A20}$$

is continuous, non-negative and satisfies $\gamma(r,0) = 0$ for all $r \geq 0$. Consequently, by virtue of Fact V in [13], there exist functions $\tilde{a} \in K_\infty$, $\tilde{\beta} \in K^+$ being increasing, such that $\gamma(r,s) \leq \tilde{a}\left(\tilde{\beta}(r)s\right)$ for all $r,s \geq 0$. Definition (A20), in conjunction with the previous inequality, implies:

$$|\Psi(t,x) - \Psi(t,y)| \leq \tilde{a}\left(\tilde{\beta}(|(t,y)|)|x-y|\right), \ \forall (t,x,y) \in \Re^+ \times \Re^m \times \Re^m \tag{A21}$$

Since $\tilde{\beta} \in K^+$ is increasing, we have $\tilde{\beta}(|(t,y)|) \leq \tilde{\beta}(2t) + \tilde{\beta}(2|y|)$ and using the properties of $K_\infty$ functions we obtain from inequality (A21):

$$|\Psi(t,x) - \Psi(t,y)| \leq \tilde{a}\left(2\tilde{\beta}(2t)|x-y|\right) + \tilde{a}\left(2\tilde{\beta}(2|y|)|x-y|\right), \ \forall (t,x,y) \in \Re^+ \times \Re^m \times \Re^m \tag{A22}$$

Define $R := 2\tilde{\beta}(0)$, $a_2(s) := s + 2\tilde{\beta}(2s) - R$ and notice $a_2 \in K_\infty$. Using again the properties of $K_\infty$ functions we obtain from inequality (A22):

$$|\Psi(t,x) - \Psi(t,y)| \leq \tilde{a}\left(2\tilde{\beta}(2t)|x-y|\right) + \tilde{a}\left(2R|x-y|\right) + \tilde{a}\left(2a_2(|y|)|x-y|\right), \ \forall (t,x,y) \in \Re^+ \times \Re^m \times \Re^m \tag{A23}$$

Inequality (2.16) is directly implied by inequality (A23) with $\beta(t) := \tilde{\beta}(2t) + R$, $a_1(s) := 2\tilde{a}(2s)$. ◁